# From Bombieri's Mean Value Theorem to the Riemann Hypothesis

Fu-Gao SONG

College of Electronic Science & Technology, Shenzhen University, China

Shenzhen Key Laboratory of Micro-nano Photonic Information Technology, China

**Abstract.** From Bombieri's mean value theorem one can deduce the prime number theorem $\pi(x) = \text{Li}(x) + O(x^{1/2} \ln^{15} x)$, which is equivalent to the Riemann hypothesis, and the least prime $P(q)$ satisfying $P(q) = O\{[\varphi(q)]^2 [\ln \varphi(q)]^{32}\}$ in arithmetic progressions with common difference $q$, where $\varphi(q)$ is the Euler function.

## 1. Introduction

In 1936, Page-Siegel and Walfisz proved the following theorem about the error term in the prime number theorem for arithmetic progressions (see also [1]):

*Suppose that $\varepsilon > 0$ is a real number and $(l, k) = 1$. Then the following estimate holds when $x \geq \exp(k^\varepsilon)$:*

$$\pi(x; l, k) = \frac{\text{Li}(x)}{\varphi(k)} + O\left( \frac{x}{\varphi(k)} \exp\left(-c(\varepsilon) \ln^{1/2} x\right) \right). \tag{1.1}$$

In 1908, Landau proved the following theorem concerning the error term of the prime number theorem, without using the theory of entire functions:

$$\pi(x) = \mathrm{Li}(x) + O\left(x\exp[-c(\ln x)^{1/2}]\right). \quad (1.2)$$

So far, the best results for the error term in the prime number theorem were obtained by Виноградов (1958) [2] and Коробов (1958) [3] using estimates of the trigonometric sum. They showed that

$$\pi(x) = \mathrm{Li}(x) + O\left(x\exp[-c(\ln x)^{3/5}(\ln\ln x)^{-1/5}]\right). \quad (1.3)$$

Many mathematicians estimate the error term in the prime number theorem by using the primary method of Selberg. There are many important results in this area of research. For example, the result that

$$\pi(x) = \mathrm{Li}(x) + O\left(\frac{x}{\ln^A x}\right), \quad (1.4)$$

was obtained by Bombieri (1962) [4] and Wirsing (1964) [5], where $A$ is any positive number. The result that

$$\pi(x) = \mathrm{Li}(x) + O\left(x\exp[-c(\ln x)^{1/7}(\ln\ln x)^{-2}]\right) \quad (1.5)$$

was obtained by Diamond and Steinig (1970) [6]; whereas

$$\pi(x) = \mathrm{Li}(x) + O\left(x\exp[-c(\ln x)^{1/6}(\ln\ln x)^{-3}]\right) \quad (1.6)$$

was obtained by Лаврик and Собиров (1973) [7].

On the other hand, von Koch [8] proved in 1901 that the following statement of the prime number theorem holds if the Riemann Hypothesis (RH) is true:

$$\pi(x) = \mathrm{Li}(x) + O\left(x^{1/2}\ln x\right); \quad (1.7)$$

in other words, if RH is true, then the prime number theorem would have the precision shown in (1.7). While according to the Generalized Riemann Hypothesis (GRH), the prime number theorem would take

$$\pi(x) = \mathrm{Li}(x) + O\left(x^{1/2}\ln^2 x\right). \quad (1.8)$$

However, the following statement of the prime number theorem is in fact equivalent to the Riemann hypothesis:

**Theorem 1.1.** *The Riemann hypothesis is equivalent to the following statement of*

*the prime number theorem*:

$$\pi(x) = \text{Li}(x) + O(x^{1/2} \ln^a x), \tag{1.9}$$

*where a > 0 is a real number.*

This theorem provides a way for proving the Riemann hypothesis.

Up to now, it has not been proved whether the discrepancy of $|\text{Li}(x) - \pi(x)|$ is less than $O(x^c)$ with $c < 1$ or not by using methods in the analytic theory of numbers.

This paper shows that a new prime number theorem that is equivalent to the Riemann Hypothesis can be obtained by using Bombieri's mean value theorem which was proved in 1965 [9].

## 2. The prime number theorem

In 1965, Bombieri proved the following theorem [9]:

**Theorem 2.1. (Bombieri's mean value theorem)** *Suppose that $A \in (0, \infty)$ and $B = A + 15$ are real numbers. The following estimate of the remainder term holds*:

$$\sum_{n \leq x^{1/2} \ln^{-B} x} \max_{y \leq x} \max_{(t,n)=1} \left| \pi(y;t,n) - \frac{\text{Li}(y)}{\varphi(n)} \right| = O\left( \frac{x}{\ln^A x} \right). \tag{2.1}$$

From Bombieri's mean value theorem one can prove the following theorem:

**Theorem 2.2. (Bombieri's mean value theorem)** *Suppose that q is a finite integer. The following estimate of the remainder term holds*:

$$\sum_{n \leq q} \max_{y \leq x} \max_{(t,n)=1} \left| \pi(y;t,n) - \frac{\text{Li}(y)}{\varphi(n)} \right| = O(x^{1/2} \ln^{15} x). \tag{2.2}$$

*Proof.* For any fixed finite integer $q$, there at least exist a real number $A$ such that $B = A + 15$ and the equation $q = [x^{1/2} \ln^{-B} x]$ holds for large enough $x$. Due to

$$O\left( \frac{x}{\ln^A x} \right) = O\left( \frac{x^{1/2}}{\ln^B x} x^{1/2} \ln^{15} x \right) = O(q x^{1/2} \ln^{15} x) = O(x^{1/2} \ln^{15} x),$$

one obtains the theorem 2.2 by substituting it into (2.1).

From Theorem 2.2 one can prove the following theorem:

**Theorem 2.3.** *Suppose that $q$ and $t$ are finite integers satisfying $(t, q) = 1$. The following estimate of the remainder term holds*:

$$\pi(x; t, q) = \frac{\text{Li}(x)}{\varphi(q)} + O(x^{1/2} \ln^{15} x). \qquad (2.3)$$

*Proof.* For any fixed finite integer $q$, Theorem 2.2 gives

$$\sum_{n \leq q-1} \max_{y \leq x} \max_{(t,n)=1} \left| \pi(y; t, n) - \frac{\text{Li}(y)}{\varphi(n)} \right| = O(x^{1/2} \ln^{15} x),$$

$$\sum_{n \leq q} \max_{y \leq x} \max_{(t,n)=1} \left| \pi(y; t, n) - \frac{\text{Li}(y)}{\varphi(n)} \right| = O(x^{1/2} \ln^{15} x).$$

Comparing above two equations one obtains the equation (2.3).

Let $q = 1$ or $q = 2$ in Theorem 2.3, one obtains the prime number theorem:

**Theorem 2.4. (The prime number theorem)** *The following estimate holds*:

$$\pi(x) = \text{Li}(x) + O(x^{1/2} \ln^{15} x). \qquad (2.4)$$

It is obvious, from Theorem 1.1, that Theorem 2.4 is in effect equivalent to the Riemann hypothesis. Therefore, Riemann hypothesis is true.

## 3. The upper bound of the least prime in arithmetic progressions

Suppose that $q$ is a finite integer, introduce the congruent integer set $\mathscr{C}(q, t)$:

$$\mathscr{C}(q, t) = \{x \mid x > 0, x \equiv t \pmod{q}, (t, q) = 1\}.$$

Let $P(q, t)$ denote the least prime in $\mathscr{C}(q, t)$ and $P(q)$ denote the maximum among all $P(q, t)$ with the same common difference $q$. One has the following theorem:

**Theorem 3.1.** *The least prime $P(q)$ in the set $\mathscr{C}(q, t)$ satisfies*:

$$P(q) = O\big([\varphi(q)]^2 \ln^{32} \varphi(q)\big). \qquad (3.1)$$

*Proof.* For any fixed $x$, there must be a finite real number $C_{qtx}$ such that equation (2.3) can be written into the following form:

$$\pi(x;t,q) = \frac{\text{Li}(x)}{\varphi(q)} + C_{qtx} x^{1/2} \ln^{15} x.$$

Let $x = C[\varphi(q)]^2 \ln^{32} \varphi(q)$ in above, where $C > 0$ is an undetermined constant. One obtains

$$\frac{\text{Li}(x)}{\varphi(q)} \sim \frac{x}{\varphi(q) \ln x} \sim \frac{C}{2} \varphi(q) \ln^{31} \varphi(q),$$

$$C_{qtx} x^{1/2} \ln^{15} x \sim 2^{15} \sqrt{C} C_{qtx} \varphi(q) \ln^{31} \varphi(q).$$

Therefore, there must be

$$\pi(x;t,q) = \frac{\text{Li}(x)}{\varphi(q)} + C_{qtx} x^{1/2} \ln^{15} x \geq \frac{\text{Li}(x)}{\varphi(q)} - |C_{qtx}| x^{1/2} \ln^{15} x > 0$$

when $x = C[\varphi(q)]^2 \ln^{32} \varphi(q)$ and $C > 2^{32} C_{qtx}^2$. The theorem follows.